\documentclass[a4paper,10pt]{scrartcl}
\usepackage[top=1in, bottom=1in, left=1in, right=1in]{geometry}
\usepackage{amssymb}
\usepackage{enumerate}
\usepackage{amsthm}
\usepackage[fleqn]{amsmath}
\usepackage[]{graphicx}
\usepackage{mathrsfs}

\newtheorem{Theorem}{Theorem}[section]
\newtheorem{Definition}[Theorem]{Definition}
\newtheorem{Axiom}[Theorem]{Axiom}
\newtheorem{Example}[Theorem]{Example}
\newtheorem{Presupposition}[Theorem]{Presupposition}

\begin{document}

\title{The importance of developing a foundation for naive category theory}
\date{}
\author{\small{Marcoen J.T.F. Cabbolet---marcoen.cabbolet@vub.ac.be}\\
        \small{Center for Logic and Philosophy of Science, Free University of Brussels}
        %\small{Marcoen.Cabbolet@vub.ac.be}
        }

\maketitle

\begin{abstract} {\footnotesize
\noindent \textbf{Abstract}---Recently Feferman (Rev. Symb. Logic \textbf{6}: 6-15, 2013) has outlined a program for the development of a foundation for naive category theory. While Ernst (\emph{ibid}. \textbf{8}: 306-327, 2015) has shown that the resulting axiomatic system is still inconsistent, the purpose of this note is to show that nevertheless \emph{some} foundation has to be developed before naive category theory can replace axiomatic set theory as a foundational theory for mathematics. It is argued that in naive category theory currently a `cookbook recipe' is used for constructing categories, and it is explicitly shown with a formalized argument that this ``foundationless'' naive category theory therefore contains a paradox similar to the Russell paradox of naive set theory.
}
\end{abstract}

\section{Introduction}
Very generally speaking, \emph{category theory} is a branch of mathematics that deals with the formalization of mathematical structures in terms of \emph{categories}. The basic notion of a category was originally introduced by Eilenberg and Maclane with the following axiomatic definition \cite{Eilenberg}:
\begin{Definition}\label{def:1}\rm A \textbf{category} $\mathscr{C}$ comprises:
\begin{enumerate}[(i)]
\item a collection of things called $\mathscr{C}$-\emph{objects}
\item a collection of things called $\mathscr{C}$-\emph{arrows}
\item for each $\mathscr{C}$-arrow $f$ a $\mathscr{C}$-object dom $f$ (the \emph{domain} of $f$) and a $\mathscr{C}$-object cod $f$ (the \emph{codomain} of $f$)
\item for each pair $\langle f, g \rangle$ of $\mathscr{C}$-arrows with dom $g$ = cod $f$, a \emph{composite} $g \circ f$ of $f$ and $g$, such that the associative law $h \circ (g \circ f) = (h \circ g) \circ f$ holds for any matching $f$, $g$, and $h$
\item for each $\mathscr{C}$-object $b$, an identity arrow $1_b$ such that the identity laws $1_b \circ f = f$ and $g \circ 1_b = g$ hold for any $\mathscr{C}$-arrow $f$ with codomain $b$ and for any $\mathscr{C}$-arrow $g$ with domain $b$. \hfill $\Box$
\end{enumerate}
\end{Definition}
\noindent To date, this is the accepted definition of a category. That being said, \emph{at least} three different positions can be discerned on the status of category theory. The first of these, championed e.g. by Mayberry \cite{Mayberry}, is that category theory needs set theory as a foundation, so that category theory itself cannot provide a foundation of mathematics. From this point of view one has to use a Bourbaki-like formalization of Def. \ref{def:1} that uses the notion of a set, such as given e.g. in \cite{Bilinsky}: the collections in the clauses (i) and (ii) of Def.~\ref{def:1} are then `sets' in the sense of ZF. The second position, originally published by Lawvere \cite{Lawvere1966}, is that category theory \emph{does} provide a foundation for mathematics. From that point of view, one is interested in categories that are not limited by the notion of a set: this is \emph{naive} or \emph{unlimited} category theory. Category theory is then considered to provide a \emph{better} foundation for mathematics than set theory, because a foundation based on categories does a better job in expressing what mathematics is about: it expresses ``the conviction that the relevant properties of mathematical objects are those which can be stated in terms of their abstract structure rather than in terms of the elements which the objects were thought to be made of'' \cite{Lawvere1966}. The third position, proposed by Landry \cite{Landry}, is that category theory provides the language of mathematics: from that point of view, category theory doesn't need set theory as a foundation, yet neither does it provide a foundation itself. The second position of these, the one espoused by Lawvere, is the topic of this paper: this is, thus, about category theory providing a foundation for mathematics.

Currently, there is no such thing as a foundation of naive category theory. Feferman, however, has outlined a program for developing one \cite{Feferman}. In short, to obtain a formal axiomatic system S that can serve as a foundation of naive category theory the following requirements should be met:
\begin{enumerate}[(R1)]
\item Form the category of all structures of a given kind, e.g. the category \emph{\textbf{Grp}} of all groups, \emph{\textbf{Top}} of all topological space, and \emph{\textbf{Cat}} of all categories.
\item Form the category $B^A$ of all functors from $A$ to $B$, where $A$, $B$ are any two categories.
\item Establish the existence of the natural numbers $N$, and carry out familiar operations on objects $a, b,\ldots$ and collections $A, B, \ldots$, including the formation of $\{a, b\}$,$(a, b)$, $A \cup B$, $A \cap B$, $A - B$, $A \times B$, $B^A$, $\bigcup A$, $\bigcap A$, $\Pi B_x [x \in A]$, etc.
\item Establish the consistency of S relative to a currently accepted system of set theory.  \hfill $\Box$
\end{enumerate}

\noindent Ernst has shown that naive category theory still is inconsistent \emph{even when} the above requirements (R1-R3) are met \cite{Ernst}: the proof is similar to the inconsistency proof following the assumption of a set of all sets $U$ (i.e. on the one hand, $U$ contains its own powerset $P(U)$ so a surjection $f: U \rightarrow P(U)$ exists; on the other hand, Cantor's theorem states that the surjection from $U$ to its powerset $P(U)$ doesn't exist). Nevertheless, the belief is that some suitable version of category theory meeting most of the above criteria is consistent (Feferman, 2015, personal communication).

The purpose of this paper is to show that one or the other foundation \textbf{necessarily} has to be in place---be it as a result of carrying out Feferman's program, be it as a result of carrying out another program---for naive category theory to be acceptable as a foundation for mathematics that replaces axiomatic set theory. The point is namely that in the present situation, in which a foundation is absent, naive category theory has to rely on a `cookbook recipe' for constructing categories: this will be formalized in the next paragraph. The final paragraph proves, with a formalized argumentation based on the Russell paradox of naive set theory, that logical inconsistency is then inevitable. So, although the Russell paradox of naive set theory is not new, this paper shows explicitly that a similar paradox is contained in naive category theory: the current ``foundationless'' naive category theory is thus worse off than axiomatic set theory as a foundation for mathematics.

\section{A formalization of the `cookbook recipe'}

We know what categories \emph{are}: if a thing $\mathscr{C}$ satisfies the axiomatic definition of a category, then $\mathscr{C}$ \emph{is} a category---at least, that's the currently prevailing view. But it remains relevant to ask how they are \emph{constructed}. Certainly, categories are not constructed the way sets are constructed in ZF: there are no constructive axioms\footnote{A \textbf{constructive axiom} is an axiom that, when certain things are given (e.g. one or two things or a thing and a predicate), states the existence of a uniquely determined other thing \cite{Bernays}.} that, when the existence of some categories is given, states the existence of a new uniquely determined category---precisely that is the aim of Feferman's program. Therefore, in naive category theory it is tacitly assumed that objects that are constructed in accordance with the axiomatic definition of a category indeed exist. In other words, it is tacitly assumed that categories can be constructed by following what can be called a \emph{cookbook recipe}. More precisely, by lack of constructive axioms the following is currently presupposed in naive category theory:
\begin{Presupposition}\label{ps:1}\rm There exists an individual category $\mathscr{C}$ for every thing constructed as follows:
\begin{enumerate}[(i)]
\item declare that there is a thing with the name $\mathscr{C}$;
\item endow the thing $\mathscr{C}$ with at least one object;
\item endow the thing $\mathscr{C}$ with at least one arrow;
\item ensure that all clauses of Def.~\ref{def:1} are satisfied;
%\item declare that there thus \emph{exists} a category $\mathscr{C}$.
\hfill $\Box$
\end{enumerate}
\end{Presupposition}
\noindent The introduction of the category \textbf{1} in Goldblatt's Topoi illustrates this ontological presupposition:
\begin{quote}
``The category \textbf{1} is to have only one object, and only one arrow. ... Suppose we call the object \textbf{a} and the arrow \emph{f}. Then we must put dom \emph{f} = cod \emph{f} = \textbf{a}, as \textbf{a} is the only object. Since \emph{f} is the only arrow, we have to take it as the identity arrow on \textbf{a}, i.e. we put $1_\textbf{a} = f$. ... This gives the identity law ... and the associative law holds ... \textbf{Thus we have a category}'' (emphasis added) \cite{Goldblatt}
\end{quote}
Other categories are endowed with far more complex collections of objects and arrows, but the principle of construction is the same as in Ps.~\ref{ps:1}: there are no proofs of existence, the categories are constructed as above and then further applied.

Lawvere already noted that ``these axioms [of Def.~\ref{def:1}] are obviously elementary (first-order)''\cite{Lawvere}, and with that in mind we intuitively formalize Def. \ref{def:1} in a first-order language \emph{without} using the notion of a set. The language then has a unary predicate $C$, with $Cx$ meaning $x$ is a \textbf{category}; a binary predicate $O$, with $O(x, y)$ meaning $y$ is an \textbf{object} of $x$; a binary predicate $A$, with $A(x, y)$ meaning $y$ is an \textbf{arrow} of $x$; unary function symbols `dom' and 'cod', with dom $x$ and cod $x$ denoting the \textbf{domain} and \textbf{codomain} of $x$; a binary function symbol $\circ$, with $g \circ f$ denoting the \textbf{composite} of $f$ and $g$; and a unary function symbol $1_{(.)}$, with $1_p$ denoting the \textbf{identity arrow} of the thing $p$. Def. \ref{def:1} can then be formalized as follows:
\begin{Definition}\label{def:2}\rm
$\forall x (Cx \Leftrightarrow A_1 \wedge A_2 \wedge A_3 \wedge A_4 \wedge A_5)$ with\footnote{The formulation of $A_4$ and $A_5$, which formalize clauses (iv) and (v) of Def. \ref{def:1}, is left as an exercise for the reader.}
\begin{enumerate}[(i)]
\item $A_1 := \exists y( O(x, y))$
\item $A_2 := \exists z (A(x, z))$
\item $A_3 := \forall f(A(x, f) \Rightarrow \exists u\exists v(O(x, u) \wedge O(x, v) \wedge u = {\rm dom}\ f \wedge v = {\rm cod}\ f)$
\item $A_4 := \ldots$
\item $A_5 := \ldots$ \hfill $\Box$
\end{enumerate}
\end{Definition}
\noindent So, Def. \ref{def:2} formalizes that a thing $x$ is a category if and only if there is at least one thing $y$ that is an object of $x$, and there is at least one thing $z$ that is an arrow of $x$, and etc. Def. \ref{def:2} is thus equivalent to Def. \ref{def:1}, without using set theory: the relations $O( . , . )$ and $A( . , . )$ are not necessarily set-theoretical $\in$-relations, and the objects and arrows are not necessarily the elements of some `set' of objects or arrows.

We now proceed by assuming that every thing $\mathscr{C}$ that can be constructed by the cookbook recipe of Ps.~\ref{ps:1} indeed \emph{exists}: we formalize Ps.~\ref{ps:1} as an existential axiom scheme, which means that we have an existential axiom for every thing $\mathscr{C}$ thus constructed. For example, we \emph{explicitly} assume that the thing \textbf{1} introduced in \cite{Goldblatt} indeed exists:
\begin{Axiom} \label{ax:1} \ \\
$\exists x(Cx \wedge x = \mathbf{1} \wedge \exists!y \exists! z(O(\mathbf{1},y)\wedge A(\mathbf{1}, z))\wedge \forall u(O(\mathbf{1},u) \Rightarrow u \neq \mathbf{1}) \wedge \forall v(A(\mathbf{1},v) \Rightarrow \exists w(O(\mathbf{1},w) \wedge v = 1_w)))$
\end{Axiom}

\noindent We also assume the existence of the things that are the categories \textbf{2}, \textbf{3}, ..., \textbf{N}, and so on, but we omit the explicit formulation of the corresponding existential axioms since it is not relevant for the main point. That being said, our formalization of the `cookbook recipe', Ps. \ref{ps:1}, is ready: we have formalized Def. \ref{def:1} in first-order language, and with a scheme of existential axioms we have formalized that every object that satisfies the axiomatic definition of a category is a category---that is, exists in the universe of categories.

\section{Proof of inconsistency}
Feferman wrote: ``There is no sensible way ... to form a category of all categories which do not belong to themselves'' \cite{Feferman}. However, now that things are formalized there \textbf{is} a sensible way to form a category of all categories which do not belong to themselves: we can now consider a category $\mathscr{R}$, which has as objects precisely those categories who are not identical to any of their objects, and which has as arrows precisely those things that are the identity arrows of the objects. The following example shows that this category can be constructed by applying the cookbook recipe of Ps.~\ref{ps:1}.
\begin{Example}\label{ex:2}\rm The category $\mathscr{R}$ is constructed as follows:
\begin{enumerate}[(i)]
\item we hereby declare that there is a thing with the name $\mathscr{R}$;
\item we hereby declare that the thing $\mathscr{R}$ is endowed with at least one object, \emph{in casu} precisely those categories that are not identical to any of their objects;
\item we hereby declare that the thing $\mathscr{R}$ is endowed with at least one arrow, \emph{in casu} precisely the identity arrows of the objects of $\mathscr{R}$;
\item all clauses of Def.~\ref{def:1} are then satisfied;
\item we hereby declare that there thus \emph{exists} a category $\mathscr{R}$. \hfill $\Box$
\end{enumerate}
\end{Example}
It can easily be proven that the thing $\mathscr{R}$ satisfies Def.~\ref{def:1}---e.g. clause (i) is satisfied because $\mathscr{R}$ has at least one object, namely the category \textbf{1}, and clause (ii) is satisfied because $\mathscr{R}$ has at least one arrow, namely the identity arrow $1_\mathbf{1}$---but we omit a formal proof: what is important is that we have this category $\mathscr{R}$ if Ps.\ref{ps:1} would hold. Proceeding, since we explicitly assume that every object that can be constructed with the cookbook recipe of Ps.~\ref{ps:1} exists, we formally assume that the thing $\mathscr{R}$ exists:
\begin{Axiom} \label{ax:2} \ \\
$\exists x(Cx \wedge x = \mathscr{R} \wedge \forall y (O(\mathscr{R}, y) \Leftrightarrow Cy \wedge \forall u(O(y, u) \Rightarrow y \neq u)) \wedge \forall z (A(\mathscr{R},z) \Leftrightarrow \exists v(O(\mathscr{R}, v) \wedge z = 1_v)))$
\end{Axiom}

\noindent The axiomatic system is then inconsistent as the following theorem shows:

\begin{Theorem}\label{th:1}\rm
$O(\mathscr{R}, \mathscr{R}) \Leftrightarrow \neg O(\mathscr{R}, \mathscr{R})$\\
\ \\
\textbf{Proof} Suppose $O(\mathscr{R}, \mathscr{R})$. But then $C\mathscr{R} \wedge \neg \forall y(O(\mathscr{R}, y) \Rightarrow \mathscr{R} \neq y)$, so from Ax.~\ref{ax:2} we then get $\neg O(\mathscr{R},\mathscr{R})$. Ergo, $O(\mathscr{R}, \mathscr{R}) \Rightarrow \neg O(\mathscr{R}, \mathscr{R})$. Next, suppose $\neg O(\mathscr{R}, \mathscr{R})$. But then $C\mathscr{R} \wedge \forall y(O(\mathscr{R}, y) \Rightarrow \mathscr{R} \neq y)$, so we get from Ax.~\ref{ax:2} that $O(\mathscr{R},\mathscr{R})$. Ergo, $\neg O(\mathscr{R}, \mathscr{R}) \Rightarrow O(\mathscr{R}, \mathscr{R})$. Q.e.d.
\end{Theorem}
%\\
%\ \\

\noindent Th.~\ref{th:1} thus demonstrates that logical inconsistency is inevitable \textbf{if} literally every thing that satisfies the currently accepted axiomatic definition of a category indeed exists as a category. Consequently, there cannot be a `cookbook recipe' that uses Def.~\ref{def:1} as an all-embracing axiomatic definition of a category for constructing categories, like: \emph{you declare that there is a thing X, you endow it with objects and arrows in accordance with Def.~\ref{def:1}, and alakazam!---there `is' your category X}. For if every category thus produced exists, logical inconsistency is inevitable.

While we do not contest any of the results obtained in category theory, we thus conclude that \emph{some} foundation for naive category theory that contains constructive axioms has to be developed, since the alternative---i.e., this ``foundationless'' naive category theory that uses the `cookbook recipe'---is inconsistent. However, \emph{even if} this Russell paradox is resolved by developing a foundation that satisfies Feferman's requirements, naive category theory \emph{still} contains an inconsistency as shown by Ernst: the question thus arises how, if at all, naive category theory can be put on a secure foundation.

\paragraph{Acknowledgement} This work has been facilitated by the Foundation Liberalitas.


\begin{thebibliography}{}

\bibitem{Eilenberg}
S. Eilenberg, S. MacLane, \emph{Trans. Amer. Math. Soc}. \textbf{58}, 231-294 (1945)

\bibitem{Mayberry}
J. Mayberry, \emph{Philosophia Mathematica} \textbf{2}(1), 16-35 (1994)

\bibitem{Bilinsky}
C. Byli\'{n}ski, \emph{Formalized Mathematics} \textbf{1}(2), 409-420 (1990)

\bibitem{Lawvere1966}
F.W. Lawvere, ``The Category of Categories as a Foundation of Mathematics'', Proceedings of the Conference on Categorical Algebra (La Jolla 1965). New York: Springer Verlag, pp. 1-20 (1966)

\bibitem{Landry}
E. Landry, \emph{Philosophy of Science} \textbf{66} (Proceedings), S14-S27 (1998)

\bibitem{Feferman}
S. Feferman, \emph{Rev. Symb. Logic} \textbf{6}(1), 6-15 (2013)

\bibitem{Ernst}
M. Ernst, \emph{Rev. Symb. Logic} \textbf{8}(2), 306-327 (2015)

\bibitem{Bernays}
P. Bernays, \emph{Axiomatic Set Theory}, Mineola: Dover Publications Inc., p. 9 (1968)

\bibitem{Goldblatt}
R. Goldblatt, \emph{Topoi: the Categorical Analysis of logic}, Mineola: Dover Publications Inc., pp. 24-27 (2006)

\bibitem{Lawvere}
F.W. Lawvere, \emph{Proc. Natl. Acad. Sci. USA} \textbf{52}, 1506-1511 (1965)

\end{thebibliography}
\end{document}